\documentclass[12pt]{amsart}
\usepackage{amsmath,amscd,amssymb,amsfonts,graphics}
\usepackage{hyperref}
\newcommand{\hl}{\hyperlink}
\newcommand{\htt}{\hypertarget}
\newcommand{\h}{\hbox}
\newcommand{\q}{\quad}

\newcommand{\bs}{\par\bigskip}
\newcommand{\ms}{\par\medskip}
\newcommand{\sk}{\par\smallskip}
\newcommand{\bsn}{\par\bigskip\noindent}
\newcommand{\msn}{\par\medskip\noindent}

\newcommand{\ges}{\geqslant}
\newcommand{\les}{\leqslant}
\newcommand{\1}{\hskip1pt}
\newcommand{\mcap}{\hbox{$\bigcap$}}
\newcommand{\mcup}{\hbox{$\bigcup$}}
\newcommand{\mopl}{\hbox{$\bigoplus$}}

\newcommand{\Cc}{{\mathcal C}}
\newcommand{\D}{{\mathcal D}}
\newcommand{\Hc}{{\mathcal H}}
\newcommand{\M}{{\mathcal M}}
\newcommand{\OO}{{\mathcal O}}
\newcommand{\ft}{\widetilde{f}}

\newcommand{\Ut}{\widetilde{U}}
\newcommand{\Xt}{\widetilde{X}}
\newcommand{\Zt}{\widetilde{Z}}
\newcommand{\Dt}{\widetilde{D}}

\newcommand{\Mt}{\widetilde{\M}}
\newcommand{\RR}{{\mathbf R}}
\newcommand{\R}{{\mathbb R}}
\newcommand{\PP}{{\mathbb P}}
\newcommand{\Q}{{\mathbb Q}}
\newcommand{\C}{{\mathbb C}}

\newcommand{\Z}{{\mathbb Z}}
\newcommand{\Gr}{{\rm Gr}}
\newcommand{\om}{\omega}
\newcommand{\Si}{\Sigma}
\newcommand{\dd}{\partial}
\newcommand{\ddd}{{\rm d}}
\newcommand{\eq}{\,{=}\,}
\newcommand{\gess}{\,{\ges}\,}

\newcommand{\mi}{\1{-}\1}
\newcommand{\mis}{\1\!{-}\!\1}
\newcommand{\pl}{\1{+}\1}
\newcommand{\bl}{\bigl}
\newcommand{\br}{\bigr}
\newcommand{\sst}{\,{\subset}\,}
\newcommand{\stm}{\,{\setminus}\,}
\newcommand{\ins}{\,{\in}\,}
\newcommand{\tos}{\,{\to}\,}
\newcommand{\ssc}{\,\raise.15ex\hbox{${\scriptstyle\circ}$}\,}
\newcommand{\ssb}{\raise.15ex\h{${\scriptscriptstyle\bullet}$}}
\newcommand{\into}{\hookrightarrow}
\newcommand{\simto}{\,\,\rlap{\hskip1.3mm\raise1.4mm\hbox{$\sim$}}\hbox{$\longrightarrow$}\,\,}
\begin{document}
\title[Proper K\"ahler morphisms]{Some remarks on decomposition theorem\\for proper K\"ahler morphisms}
\author{Morihiko Saito}
\address{RIMS Kyoto University, Kyoto 606-8502 Japan}
\begin{abstract} We explain a correct proof of the decomposition theorem for direct images of constant Hodge modules by proper K\"ahler morphisms of complex manifolds. We also give some examples showing certain difficulty in the non-constant Hodge module case.
\end{abstract}
\maketitle
\centerline{\bf Introduction}
\bsn
Decomposition theorem for direct images of {\it constant\1} Hodge modules under proper K\"ahler morphisms of complex manifolds was announced in \cite[Thm.\,0.5]{toh}. There was, however, some misstatement in the proof. In this paper we explain a correct proof of the following.
\par\htt{T1}{}\msn
{\bf Theorem 1.} {\it Let $f:X\tos Y$ be a proper K\"ahler morphism of connected complex manifolds with $\dim X\eq n$. Then the Hodge filtration $F$ on the direct image as a filtered $\D$-module $f^{\D}_*(\OO_X,F)$ is strict, the cohomological direct images
$$\Hc^j\!f_*\bl((\OO_X,F),\R_X[n]\br):=\bl(\Hc^j\!f^{\D}_*(\OO_X,F),{}^p\!R^j\!f_*\R_X[n]\br)\q(j\ins\Z)$$
are Hodge modules of weight $j\mi n$, we have the hard Lefschetz property
\htt{1}{}
$$\ell^k:\Hc^{-k}\!f_*\bl((\OO_X,F),\R_X[n]\br)\simto\Hc^k\!f_*\bl((\OO_X,F),\R_X[n]\br)(k)\q(k\,{>}\,0),
\leqno(1)$$
and the primitive part is polarized with sign given as in \cite[5.3.1 and 5.4.1]{mhp}. Here $\ell$ is the cohomology class of a relative K\"ahler class in $H^2(X,\R(1))$, and $\Gr^F_p\OO_X=0$ $(p\ne 0)$.}
\ms
This implies the decomposition theorem using \cite{De1} together with the strict support decomposition. The difficulty in the proof of Theorem~\hl{T1}{1} is that it is not easy to deduce the {\it hard Lefschetz\1} property~(\hl{1}{1}) from the one for the composition with a desingularization. We may assume that the latter is a composition of blow-ups along smooth centers, see \cite{Wl}. In the {\it constant\1} Hodge module case, this situation can be utilized quite effectively, since we can easily apply the key lemma \cite[Lem.\,5.2.15]{mhp} to the restriction and Gysin morphisms between the blow-up and the exceptional divisor at each step. It does not, however, seem easy to generalize this to the {\it non-constant\1} case, since extension classes of general Hodge modules can be rather complicated, see \hl{S2.3}{2.3} below. (This seems to imply that some argument in twistor theory does not necessarily apply to Hodge modules.)
\sk
As a corollary we get the following (extending the main theorem in \cite{FFS} to the proper K\"ahler case).
\par\htt{T2}{}\msn
{\bf Theorem 2} (see \cite[Thm\,2.5]{Fn}). {\it Let $(X,D)$ be an analytic simple normal crossing pair with $D$ reduced. Let $f\,{:}\,X\tos Y$ be a proper morphism to a complex manifold $Y$. Assume $f$ is K\"ahler on each irreducible component of $X$. Then we have the weight spectral sequence
\htt{2}{}
$${}_FE_1^{-q,i+q}=\mopl_{k+l=n+q+1}\,R\1^i\!f_*\om_{D^{[k,l]}/Y}\Longrightarrow R\1^i\!f_*\om_{X/Y}(D),
\leqno(2)$$
degenerating at $E_2$, and its $E_1$-differential $\ddd_1$ splits so that the $E_2^{-q,i+q}$ are direct factors of $E_1^{-q,i+q}$.}
\ms
In Section 1 we prove Theorem~\hl{T1}{1} with $X$ replaced by an embedded resolution shrinking $Y$ if necessary. This implies the theorem in the general case except the hard Lefschetz property and the induced polarization on the $\ell$-primitive part. In Section 2 we prove the hard Lefschetz property after reviewing some basics of Gysin and restriction morphisms associated to blow-ups of complex manifolds along smooth centers.
\ms
We thank O.~Fujino for his interest in this subject. This work is partially supported by JSPS Kakenhi 15K04816.
\bs\bs\htt{S1}{}
\vbox{\centerline{\bf 1. Proof of Theorem~\hl{T1}{1} for embedded resolutions}
\bsn
In this section we prove Theorem~\hl{T1}{1} with $X$ replaced by an embedded resolution shrinking $Y$ if necessary. This implies the theorem in the general case except the hard Lefschetz property and the induced polarization on the $\ell$-primitive part.}
\par\htt{1.1}{}\msn
{\bf 1.1.~Level $0$ Hodge modules of normal crossing type.} Let $X$ be a complex manifold of pure dimension $n$. We say that a polarizable $\R$-Hodge module $\M\eq\bl((M,F),K\br)$ with strict support $X$ (that is, $K$ is an intersection complex on $X$ with local system coefficients, see \cite{BBD}) is a {\it level $0$ Hodge module of normal crossing type\1} with singular locus ${\rm Sing}\,\M$ contained in $D\sst X$ if $D$ is a divisor with {\it simple\1} normal crossings and the restriction of $\M$ to $U:=X\stm D$ corresponds to a variation of Hodge structure of level 0 and type $(0,0)$ (in particular, $\M$ has weight $n$). The last condition means that the Hodge filtration $F$ is trivial on $U$, more precisely, $\Gr^F_pM|_U\eq0$ for $p\,{\ne}\,0$. The underlying local system $K[-n]|_U$ is an orthogonal representation (by a polarization) with semisimple local monodromies of finite orders around the irreducible components of $D$. Note that the constant Hodge module $\bl((\OO_X,F),\R_X[n]\br)$ is a level 0 Hodge module of normal crossing type with weight $n$, where $D\eq\emptyset$.
\sk
Let $\M$ be a level 0 Hodge module of normal crossing type with ${\rm Sing}\,\M\sst D$. Let $f\,{:}\,X\tos Y$ be a proper K\"ahler morphism of complex manifolds. Let $p\ins f(X)$, and $g_i$ $(i\in[1,r])$ be holomorphic functions on $Y$ (replacing $Y$ with a neighborhood of $p$ if necessary) such that
\htt{1.1.1}{}
$$\mcap_{i=1}^r\,g_i^{-1}(0)\cap f(X)=\{p\},
\leqno(1.1.1)$$
where we may have $r>\dim f(X)$. Consider the condition
\htt{1.1.2}{}
$$D\supset\mcup_{i=1}^r\,h_i^{-1}(0),
\leqno(1.1.2)$$
with $h_i:=f{}^*g_i$ (replacing $D$ and taking an embedded resolution). Here the divisors $h_i^{-1}(0)$ have simple normal crossings, since they are contained in a divisor with simple normal crossings. Let $\ell$ be the cohomology class of a relative K\"ahler form for $f$. We have the following.
\par\htt{T1.1}{}\msn
{\bf Theorem~1.1.} {\it In the above notation, the direct image as a filtered $\D$-module $f^{\D}_*(M,F)$ is strict, and the cohomological direct images $\Hc^j\!f_*\M$ $(j\ins\Z)$ are Hodge modules of weight $j{-}n$, where {\rm(\hl{1.1.2}{1.1.2})} is not assumed. If {\rm(\hl{1.1.2}{1.1.2})} is satisfied, then, replacing $Y$ with a neighborhood of $p$ if necessary, we have the hard Lefschetz property for the action of $\ell$ together with the induced polarization on the $\ell$-primitive part.}
\msn
{\it Proof.} We argue by induction on $\dim f(X)$. If $f(X)\eq\{p\}$, we may assume $Y\eq\{p\}$. Hodge theory for unitary local systems with semi-simple local monodromies of finite orders is an easy case of \cite{CKS2}, \cite{KK1}, \cite{KK2} (generalizing \cite{Zu}), and is more or less well-known, see for instance \cite[\S5]{Ti}.
\sk
Assume $\dim f(X)\,{>}\,0$. Applying \cite[Thm.\,3.3--4]{mhm}, the hypothesis that $\M$ is a level 0 Hodge module of normal crossing type is preserved by taking the $N$-primitive part of the $W$-graded pieces of nearby and vanishing cycle functors $P_N\Gr^W_k\psi_{h_1}$, $P_N\Gr^W_k\varphi_{h_1,1}$ up to Tate twists. Assuming (\hl{1.1.2}{1.1.2}), and applying \cite[Prop.3.3.17, Prop.4.2.2, and Cor.4.2.4]{mhp} and the inductive hypothesis, we then see that
$$P_N\Gr^W_k\psi_{g_1}\Hc^jf_*\M,\q P_N\Gr^W_k\varphi_{g_1,1}\Hc^jf_*\M\q(k\ins\Z)$$
are Hodge modules of weight $k$ (using the weight spectral sequence), and get the strictness of the Hodge filtration $F$ on the direct image together with the hard Lefschetz property and the induced polarization on the $(\ell,N)$-bi-primitive part, replacing $Y$ with a neighborhood of $p$ if necessary. (Note that the vanishing cycle functor is the identity for objects whose supports are contained in the zero-locus of the function.)
\sk
Here we can apply the above argument by replacing $\psi_{g_1}$, $\varphi_{g_1,1}$ with $\psi_{g_1-c}$, $\,\varphi_{g_1-c,1}$ for $0\,{<}\,|c|\,{\ll}\,1$. (Note that any intersections of irreducible components of $D$, which are not contained in $h_1^{-1}(0)$, intersects $h_1^{-1}(c)$ transversally for $0\,{<}\,|c|\,{\ll}\,1$.) This can be used to show that we have the induced polarization on the $\ell$-primitive part. Indeed, if a strict support $Z$ of a direct factor of $\Hc^j\!f_*\M$ is not contained in $g_1^{-1}(0)$ and $Z\cap g^{-1}(0)\ne\emptyset$, then we take $c$ such that $0\,{<}\,|c|\,{\ll}\,1$ and $Z^o\cap g_1^{-1}(c)\ne\emptyset$, where $Z^o$ is a smooth Zariski-open of $Z$ on which we have a shifted local system.
\sk
The above argument is, however, insufficient to prove the assertion that the cohomological direct images $\Hc^j\!f_*\M$ are Hodge modules. (The last argument on the induced polarization is valid only after showing the latter assertion.) To show it, we have to take any point $p'\ins f(X)$ together with holomorphic functions $g'_i$ $(i\in[1,r])$ satisfying (\hl{1.1.1}{1.1.1}--\hl{1.1.2}{2}) (with $p$ replaced by $p'$) after replacing $Y$ with a neighborhood of $p'$ and taking an embedded resolution, which is given by a {\it projective\1} morphism $\pi\,{:}\,\Xt\tos X$. (Here we do not have to assume (\hl{1.1.2}{1.1.2}) for $p$.) For the last morphism, we have the compatibility of decompositions for filtered $\D$-modules and $\R$-complexes (see \cite{toh}) using Deligne's {\it canonical choice\1} of decomposition depending on $\ell$, see \cite{De2}). So $\Hc^j\!f_*\M$ is a direct factor of $\Hc^j\ft_*\Mt$ with $\ft:=f\ssc\pi$, where $\Mt$ is the level 0 Hodge module of normal crossing type on $\Xt$ whose restriction to $\Xt\stm\pi^{-1}(D)$ coincides with the restriction of $\M$. Since the first two assertions are local on $Y$, Theorem~\hl{T1.1}{1.1} then follows.
\par\htt{R1.1a}{}\msn
{\bf Remark~1.1a.} In \cite{toh} an argument similar to \cite{ast} was used.
\par\htt{R1.1b}\msn
{\bf Remark~1.1b.} There is a subtle point in a proof of \cite[Thm.\,5.3.1]{mhp} and \cite[Thm.\,3.20]{mhm} by induction on the dimension of support without using \cite[3.5]{mhp}. Indeed, for the proof of \cite[Thm.\,3.20]{mhm}, we must replace $X$ with a desingularization for each choice of a function $g$, and the inductive argument cannot be carried out so simply as is imagined by some people. We have to use a rather technical assertion like \cite[Prop.\,5.3.4]{mhp} (or \cite[Thm.\,2.4]{ypg}).
\par\htt{1.2}{}\msn
{\bf 1.2.~Proof of Theorem~2.} We show that Theorem~\hl{T2}{2} follows from Theorem~\hl{T1.1}{1.1} without relying on Theorem~\hl{T1}{1}. Except the splitting of $\ddd_1$, the assertions of Theorem~\hl{T2}{2} are local on $Y$, and follow from Theorem~\hl{T1.1}{1.1}. Here we need the hard Lefschetz property for the action of a relative K\"ahler class in order to get a global splitting of $\ddd_1$ (using the associated global polarization of variation of Hodge structure).
\sk
The latter hard Lefschetz property is, however, needed only for the direct factor of the cohomological direct images with strict support $f(X)$ in the notation of Theorem~\hl{T1.1}{1.1}. Indeed, the higher direct images $R^if_*\om_{D^{[k,l]}/Y}$ are uniquely determined by their restrictions to a dense Zariski-open subset of $f(D^{[k,l]})$ over which $f|_{D^{[k,l]}}$ is smooth (using the filtration $V$ of Kashiwara and Malgrange, see \cite[(3.2.2.2)]{mhp}). This is proved in \cite{kcon} as a refinement of Koll\'ar's {\it torsion-freeness\1} theorem \cite{Ko}.
\sk
The hard Lefschetz property for {\it general fibers\1} of $f|_{D^{[k,l]}}\,{:}\,D^{[k,l]}\tos f(D^{[k,l]})$ follows from the classical Hodge theory on compact smooth K\"ahler general fibers. The weight spectral sequence (\hl{2}{2}) can be deduced from the corresponding one for Hodge modules (by restricting to $F_0$ and taking the tensor product with $\om_Y^{\vee}$, see \cite{FFS}), where the differential $\ddd_1$ {\it preserves the strict supports}. So Theorem~\hl{T2}{2} follows.
\par\htt{R1.2}{}\msn
{\bf Remark 1.2.} We say that a proper morphism of complex manifolds $f:X\to Y$ is K\"ahler, if there is a relative K\"ahler form $\xi_f$ (which is a closed real 2-form on $X$) such that locally on $Y$ there is a K\"ahler form $\xi_Y$ such that $\xi_f\1{+}\1f^*\xi_Y$ is a K\"ahler form on $X$. (Here ``locally on $Y$" can be replaced with ``locally on $X$" as in \cite[Def.\,6.1]{Ta}, \cite{Va} using semi-positivity of pullbacks of K\"ahler forms on $Y$ since $f$ is proper.) This definition seems apparently different from the one in \cite[Def.\,4.1]{Fk}. (Recall that a $C^{\infty}$ function $\phi$ on a complex manifold is {\it strictly plurisubharmonic\1} if and only if $i\dd\overline{\dd}\phi$ is K\"ahler, see for instance \cite[Rem.\,1.1]{Fk}.)
There is a problem if one wants to extend the above definition to the $X$ singular case, since the {\it cohomology class\1} of $\xi_f$ is not easy to define.
\bs\bs\htt{S2}{}
\vbox{\centerline{\bf 2. Hard Lefschetz property}
\bsn
In this section we prove the hard Lefschetz property after reviewing some basics of Gysin and restriction morphisms associated to blow-ups of complex manifolds along smooth centers.}
\par\htt{2.1}{}\msn
{\bf 2.1.~Gysin and restriction morphisms.} Let $\pi:\Xt\to X$ be the blow-up of a complex manifold $X$ along a closed submanifold $Z\sst X$ with codimension $d$. Set $\Zt:=\pi^{-1}(Z)$ with $i_{\Zt}\,{:}\,\Zt\,{\into}\,\Xt$ the inclusion. This is a $\PP^{d-1}$-bundle over $Z$. We have the restriction and Gysin morphisms
$$\rho:\Q_{\Xt}\to\Q_{\Zt},\q\gamma:\Q_{\Zt}\to\Q_{\Xt}(1)[2]\q\h{in}\,\,\,\,D^b_c(\Xt,\Q),$$
where the direct image by the closed immersion $i_{\Zt}$ is omitted to simplify the notation. Set $n\,{:=}\,\dim X$, $m\,{:=}\,\dim Z\,(=n\mi d)$.
\sk
These morphisms can be embedded into the distinguished triangles
\htt{2.1.1}{}
$$\aligned(j_{\Ut})_!\Q_{\Ut}&\to\Q_{\Xt}\buildrel{\rho}\over\to\Q_{\Zt}\buildrel{+1}\over\to,\\ \Q_{\Zt}(-1)[-2]&\buildrel{\gamma}\over\to\Q_{\Xt}\to\RR(j_{\Ut})_*\Q_{\Ut}\buildrel{+1}\over\to,\endaligned
\leqno(2.1.1)$$
which are dual of each other, where $\Ut:=\Xt\stm\Zt\,(\eq X\stm Z)$ with $j_{\Ut}\,{:}\,\Ut\,{\into}\,\Xt$ the inclusion.
\sk
It is easy to verify that the compositions $\gamma\ssc\rho$ and $\rho\ssc\gamma$ are respectively identified with the actions of the cohomology class of $\Zt\sst\Xt$ and its restriction to $\Zt$, using the canonical isomorphisms
\htt{2.1.2}{}
$$H^2(\Xt,\Q)(1)={\rm Hom}\bl(\Q_{\Xt},\Q_{\Xt}(1)[2]\br)\,\,\h{(and similarly for $\Zt$)}.
\leqno(2.1.2)$$
\sk
By the decomposition theorem for projective morphisms \cite{mhp}, we have the isomorphisms
\htt{2.1.3}{}
$$\aligned\RR\pi_*\Q_{\Xt}[n]&\,\cong\,\Q_X[n]\,\oplus\,\mopl_{i=0}^{d-2}\,\bl(\Q_Z[m]\br)(-i{-}1)[d{-}2{-}2i],\\ \RR\pi_*\Q_{\Zt}[n{-}1]&\,\cong\,\mopl_{i=0}^{d-1}\,\bl(\Q_Z[m]\br)(-i)[d{-}1{-}2i].\endaligned
\leqno(2.1.3)$$
Here we can take Deligne's canonical choice of decomposition using the relative ample divisor $-\Zt$, see \cite{De2}. (It is not clear whether these direct sum decompositions are compatible with the Gysin and restriction morphisms.) Fixing the decomposition, we can define the decreasing filtration $G^{\ssb}$ on $\RR\pi_*\Q_{\Zt}[n{-}1]$ so that
\htt{2.1.4}{}
$$G^k\1\RR\pi_*\Q_{\Zt}[n{-}1]\,\cong\,\mopl_{i=k}^{d-1}\,\bl(\Q_Z[m]\br)(-i)[d{-}1{-}2i].
\leqno(2.1.4)$$
This is opposite to a well-defined filtration $G'_k$ ($k\ins\Z)$ defined by the truncations $\tau_{\les 1-n+2k}$ (or $^{\bf p}\tau_{\les 1-d+2k}$, see \cite{BBD}).
The first decomposition of (\hl{2.1.3}{2.1.3}) can be induced by using the restriction morphism for $\Xt\tos X$ and the Gysin morphism for $\Zt\,{\into}\,\Xt$\,:
\htt{2.1.5}{}
$$\aligned\Q_X[n]&\to\RR\pi_*\Q_{\Xt}[n],\\ \bl(G'_{d-2}\1\RR(\pi_{\Zt})_*\Q_{\Zt}\br)(-1)[n{-}2]&\to\RR\pi_*\RR\Gamma_{\Zt}\Q_{\Xt}[n]\to\RR\pi_*\Q_{\Xt}[n].\endaligned
\leqno(2.1.5)$$
\par\htt{2.2}{}\msn
{\bf 2.2.~Proof of the hard Lefschetz property.} Let $f\,{:}\,X\tos Y$ be a proper K\"ahler morphism of connected complex manifolds. Set
$$\M_X:=\bl((\OO_X,F),\R_X[\dim X]\br)\,\,\,\h{(similarly for $\M_{\Xt},\M_D,\M_{\Dt})$.}$$
Put $f_Z:=f|_Z$, $\ft_{\Zt}:=\ft|_{\Zt}$ with the notation of \hl{2.1}{2.1}. By Theorem~\hl{T1.1}{1.1}, the cohomological direct images
$$\Hc^j\!f_*\M_X,\,\,\Hc^j\ft_*\M_{\Xt},\,\,\,\Hc^j(f_Z)_*\M_Z,\,\,\,\Hc^j(\ft_{\Zt})_*\M_{\Zt}$$
are locally polarizable Hodge modules on $Y$ of weight $n{+}j$, $n{+}j$, $m{+}j$, $n{-}1{+}j$ respectively. Using \cite{De2} as is explained in \cite{toh}, we can get the isomorphisms
\htt{2.2.1}{}
$$\aligned&\Hc^j\ft_*\M_{\Xt}\,\cong\,\Hc^j\!f_*\M_X\,\oplus\,\mopl_{i=0}^{d-2}\,\Hc^{j+d-2-2i}(f_Z)_*\M_Z(-i{-}1),\\&\Hc^j(\ft_{\Zt})_*\M_{\Zt}\,\cong\,\mopl_{i=0}^{d-1}\,\Hc^{j+d-1-2i}(f_Z)_*\M_Z(-i).\endaligned
\leqno(2.2.1)$$
These are compatible with the direct image of the decompositions (\hl{2.1.3}{2.1.3}) (using \cite{De2}).
\sk
Let $\ell$ be the pull-back to $\Xt$ of a relative K\"ahler class $\ell_X$ for $f$. Let $\ell'$ be the cohomology class of the relatively ample divisor $-\Zt$ for $\pi$. It is known that $\ell_c\,{:=}\,\ell\pl c\1\ell'$ ($0\,{<}\,c\,{\ll}\,1$) is a relative K\"ahler class for $\ft$.
\sk
For the proof of the hard Lefschetz property (\hl{1}{1}) and the induced polarization on the primitive part, it is enough to show the following.
\par\htt{T2.2}{}\msn
{\bf Theorem~2.2.} {\it Assume the hard Lefschetz property {\rm(\hl{1}{1})} and the induced polarization on the primitive part hold for $\ft$, $\ell_c$ $(0\,{<}\,c\,{\ll}\,1)$ and also for $f_Z$, $\ell_Z$, where $f_Z:=f|_Z$, $\ell_Z:=\ell_X|_Z$. Then they hold for $f$ and $\ell_X$.}
\msn
{\it Proof.} For $j\ins\Z$, set
$$\M_{\Xt}^{(j)}:=\Hc^{-j}\ft_*\M_{\Xt}(-j)\,\,\,\,\,\h{(similarly for $\M_{\Zt}^{(j)}$, etc.)}$$
By (\hl{2.2.1}{2.2.1}) we have the isomorphisms
\htt{2.2.2}{}
$$\aligned&\M_{\Xt}^{(j)}\,\cong\,\M_X^{(j)}\,\oplus\,\mopl_{i=0}^{d-2}\,M_Z^{(j-d+2+2i)}(i{-}d{+}1),\\&\M_{\Zt}^{(j)}\,\cong\,\mopl_{i=0}^{d-1}\,\M_Z^{(j-d+1+2i)}(i{-}d{+}1).\endaligned
\leqno(2.2.2)$$
Put $\M_{\Xt}^{(\ssb)}:=\mopl_{j\in\Z}\,\M_{\Xt}^{(j)}$, etc. We have the actions
$$\ell,\ell':\M_{\Xt}^{(\ssb)}\to\M_{\Xt}^{(\ssb)}(-1)[-2],$$
and similarly for $\M_{\Zt}^{(\ssb)}$.
\sk
By definition (see \cite[5.4.1]{mhp}) the polarizations of Hodge modules
\htt{2.2.3}{}
$$\aligned S_{\Xt}:\R_{\Xt}[n]&\otimes_{\R}\R_{\Xt}[n]\to\R_{\Xt}[2n]\,\bl(=a_{\Xt}^!\R(-n)\br),\\ S_{\Dt}:\R_{\Dt}[n{-}1]&\otimes_{\R}\R_{\Dt}[n{-}1]\to\R_{\Dt}[2n{-}2]\,\bl(=a_{\Dt}^!\R(1{-}n)\br),\endaligned
\leqno(2.2.3)$$
are defined by $(-1)^{n(n-1)/2}$ and $(-1)^{(n-1)(n-2)/2}$ respectively, where $a_{\Xt}:\Xt\tos{\rm pt}$ denotes the constant morphism.
\sk
Let $K_{\Xt}^{(j)}$ be the underlying $\R$-complex of $\M_{\Xt}^{(j)}$, and similarly for $K_{\Dt}^{(j)}$ ($j\ins\Z$). We have the induced pairings
\htt{2.2.4}{}
$$\aligned S_{\Xt}^{(j)}:K_{\Xt}^{(j)}\otimes_{\R}K_{\Xt}^{(-j)}&\to a_Y^!\R(-n),\\ S_{\Dt}^{(j)}:K_{\Dt}^{(j)}\otimes_{\R}K_{\Dt}^{(-j)}&\to a_Y^!\R(1{-}n),\endaligned
\leqno(2.2.4)$$
{\it where the sign is twisted by\1} $(-1)^{j(j-1)}$. This is closely related to the sign in \cite[(5.3.1.3)]{mhp}. Because of this change of sign, we can get the equalities
\htt{2.2.5}{}
$$S^{(j)}_{\Xt}\ssc({\rm id}{\otimes}\ell)=-S^{(j-2)}_{\Xt}\ssc(\ell{\otimes}{\rm id})\q\h{(similarly for $\ell'$)},
\leqno(2.2.5)$$
as morphisms from $K_{\Xt}^{(j)}{\otimes}K_{\Xt}^{(2-j)}$ to $a_Y^!\R(-n{-}1)$ in the derived category (similarly with $S^{(\ssb)}_{\Xt}$ replaced by $S^{(\ssb)}_{\Dt}$), see also \cite{De3}, \cite[3.4--5]{GN} for related arguments.
\sk
Let $L$ be the monodromy filtration for the action of $\ell$ on $\M_{\Xt}^{(\ssb)}$, $\M_{\Zt}^{(\ssb)}$. (This commutes with $\gamma,\rho$, since this is the action of a cohomology class of $\Xt$.) By \cite{CK}, \cite{CKS1}, we can apply \cite[Lem.\,5.2.15]{mhp} to fibers of generic variations of Hodge structure for each strict support of $\M_{\Xt}^{(\ssb)}$, $\M_{\Zt}^{(\ssb)}$ {\it using the $\ell$-primitive decomposition as graded $\R[\ell']$-modules.} We then get the decomposition
\htt{2.2.6}{}
$$\Gr^L_{\ssb}\M_{\Xt}^{(\ssb)}={\rm Ker}\,\Gr^L_{\ssb}\rho\oplus{\rm Im}\,\Gr^L_{\ssb}\gamma,
\leqno(2.2.6)$$
where
\htt{2.2.7}{}
$${\rm Im}\,\Gr^L_{\ssb}\rho=\Gr^L_{\ssb}G^1\M_{\Xt}^{(\ssb)},\q{\rm Ker}\,\Gr^L_{\ssb}\gamma=\Gr^L_{\ssb}G^{d-1}\M_{\Xt}^{(\ssb)},
\leqno(2.2.7)$$
with $G$ as in (\hl{2.1.4}{2.1.4}). (Note that these never hold without taking $\Gr^L_{\ssb}$.) Here it is not necessarily easy to verify the {\it sign\1} in general, see \cite[Rem.\,(iii) in 2.6]{ypg}. If the sign is wrong, the lemma would imply a decomposition of $\Gr^L_{\ssb}\M_{\Zt}^{(\ssb)}$ instead of $\Gr^L_{\ssb}\M_{\Xt}^{(\ssb)}$. However, one can verify that this is not possible in this {\it constant sheaf\1} case. Indeed, the sign should be essentially {\it local\1} on $\Xt$ using \cite[Thm.\,5.3.1]{mhp}, and should depend only on the dimension of $\Xt$ in the constant sheaf case, see Remark~\hl{R2.2}{2.2} below
\sk
By hypothesis, the hard Lefschetz property holds for $f_Z$, $\ell_Z$. So the filtration $L$ on $\M_{\Zt}^{(\ssb)}$ is induced via the second isomorphism of (\hl{2.2.2}{2.2.2}) by the filtration $L$ on $\M_Z^{(\ssb)}$ such that
\htt{2.2.8}{}
$$L_k\M_Z^{(\ssb)}=\mopl_{j\les m+k}\,\M_Z^{(j)},
\leqno(2.2.8)$$
since $\ell$ comes from $\ell_X$. It is then easy to verify that the action of $\ell'$ on $\Gr^L_k\M_{\Zt}^{(\ssb)}$ is essentially given by the identity (or zero) via the second isomorphism of (\hl{2.2.2}{2.2.2}) as is expected.
\sk
Since $\ell$ is the pullback of $\ell_X$, its action is compatible with the decomposition in the first isomorphism of (\hl{2.2.2}{2.2.2}), and we can apply a similar argument to the second direct factor. We then see that ${\rm Im}\,\Gr^L_{\ssb}\gamma$ coincides with the second direct factor using the decomposition (\hl{2.2.6}{2.2.6}). Indeed, ${\rm Im}\,\Gr^L_{\ssb}\gamma$ contains the second direct factor by an argument around (\hl{2.1.5}{2.1.5}), and ${\rm Im}\,\Gr^L_{\ssb}\gamma$ is injectively sent to $\M_{\Zt}^{(\ssb)}$ by $\Gr^L_{\ssb}\rho$ as a consequence of (\hl{2.2.6}{2.2.6}) (and $\rho\ssc\gamma$ coincides with the action of $\ell'$ up to a sign). We then get the isomorphism as bi-graded Hodge modules
\htt{2.2.9}{}
$$\Gr^L_{\ssb}\M_X^{(\ssb)}\,\cong\,{\rm Ker}\,\Gr^L_{\ssb}\rho\,\bl(\subset\,\Gr^L_{\ssb}\M_{\Xt}^{(\ssb)}\br).
\leqno(2.2.9)$$
\sk
For the proof of the hard Lefschetz property for $f$, $\ell_X$, it is enough to show
\htt{2.2.10}{}
$$\Gr^L_k\M_X^{(j)}=0\,\,\,(k\ne j\pl n).
\leqno(2.2.10)$$
By \cite{CK}, \cite{CKS1} together with the hypothesis, the relative monodromy filtration for $\ell',L$ coincides with the monodromy filtration for $\ell_c\,{:=}\,\ell\pl c\1\ell'$ ($0\,{<}\,c\,{\ll}\,1$), and the hard Lefschetz property holds for this so that the monodromy filtration for $\ell_c$ coincides with the filtration by degree up to shift. By definition, the relative monodromy filtration induces the monodromy filtration for $\Gr^L_k\ell'$ shifted by $k$ on each $\Gr^L_k\M_{\Xt}^{(\ssb)}$. It is easy to see that the action of $\Gr^L_{\ssb}\ell'$ on $\Gr^L_{\ssb}\M_{\Xt}^{(\ssb)}$ is compatible with the direct sum decomposition (\hl{2.2.6}{2.2.6}), since $\gamma\ssc\rho\eq{\pm}\1\ell'$. The action of $\Gr^L_{\ssb}\ell'$ vanishes on the direct factor ${\rm Ker}\,\Gr^L_{\ssb}\rho$, hence the relative monodromy filtration (or the monodromy filtration for $\ell_c$) induces a trivial filtration on each ${\rm Ker}\,\Gr^L_k\rho$. So the condition (\hl{2.2.10}{2.2.10}) is satisfied for ${\rm Ker}\,\Gr^L_{\ssb}\rho$, and it holds for $\Gr^L_{\ssb}\M_X^{(\ssb)}$ using the bi-graded isomorphism (\hl{2.2.9}{2.2.9}).
\sk
It remains to show the assertion on the induced polarization on the $\ell$-primitive part.
For this we can apply a well-known assertion saying that a continuous family of {\it non-degenerate\1} hermitian matrices parametrized by a connected interval in $\R$ is everywhere positive-definite if so is at one point. This finishes the proof of Theorem~\hl{T2.2}{2.2}.
\par\htt{R2.2}\msn
{\bf Remark~2.2.} The sign problem in the decomposition (\hl{2.2.6}{2.2.6}) is not difficult in the {\it constant\1} Hodge module case. Indeed, by \cite[Thm.\,5.3.1]{mhp}, the problem is the signs related to the polarizations of Hodge modules (\hl{2.2.3}{2.2.3}) and to the restriction and Gysin morphisms
$$\rho\in{\rm Ext}^1\bl(\R_{\Xt}[n],\R_{\Dt}[n{-}1]\br),\q\gamma\in{\rm Ext}^1\bl(\R_{\Dt}(-1)[n{-}1],\R_{\Xt}[n]\br).$$
(Note that their extension classes correspond respectively to $j_!j^*\R_{\Xt}[n]$ and $\RR j_*j^*\R_{\Xt}[n]$ with $j\,{:}\,\Xt\stm\Dt\,{\into}\,\Xt$ the inclusion.) More precisely, the problem is whether we have
\htt{2.2.11}{}
$$\aligned S_{\Xt}\ssc(\gamma\otimes{\rm id})&=\gamma\ssc S_{\Dt}\ssc({\rm id}\otimes\rho)\\ \h{or}\q\gamma\ssc S_{\Dt}\ssc(\rho\otimes{\rm id})&=S_{\Xt}\ssc({\rm id}\otimes\gamma),\endaligned
\leqno(2.2.11)$$
since $S_{\Xt}$, $S_{\Dt}$ are $(-1)^n$ and $(-1)^{n-1}$-{\it symmetric.} These morphisms are defined over $\Z$, and the Gysin morphism $\gamma$ is a kind of trace morphism. In the case the first (resp.\ second) equality holds, we say that there is a $\M_{\Xt}$ (resp.\ $\M_{\Dt}$)-{\it decomposable relation.} (This is closely related to \cite[Lem.\,5.2.15]{mhp}.)
\sk
The above relation is stable by the direct image under a projective morphism using \cite[Thm.\,5.3.1]{mhp} (restricting to the primitive part and using the action of a relative ample class). The problem is then {\it local\1} on $\Xt$ by the adjunction relation between the direct image and pullback under the inclusion $\Dt\,{\into}\,\Xt$, and is {\it independent of\1} $f$. In this constant Hodge module case, we then see that there is a $\M_{\Dt}$-decomposable relation by considering the case $\Xt$ is a smooth projective variety with $\Dt$ a sufficiently general hyperplane section as a corollary of the weak Lefschetz theorem. In the case $\Dt$ is the exceptional divisor of a blow-up along a smooth center, we have the converse relation, since $-\Dt$ is a relative ample divisor, and the minus sign here reverses the relation. Note that the sign becomes much more complicated if the (Hermitian) dual is used instead of {\it perfect pairings of constructible complexes\1} (as in the twistor case), see also \cite[Rem.\,1.1e]{FPS}.
\par\htt{S2.3}\msn
{\bf 2.3.~Generalization to the non-constant case.} A similar argument in a much more general situation seems to be employed in twistor theory. It does not seem, however, quit easy to generalize the above argument to the case of {\it non-constant\1} Hodge modules, since the extension groups of non-constant Hodge modules can be rather complicated. Here we have to construct first a generalization of restriction and Gysin morphisms.
\sk
Let $\Xt$ be a connected complex manifold with $\Dt\sst\Xt$ a divisor with simple normal crossings. Let $\M_{\Xt}$ be a pure Hodge module of weight $w$ with strict support $\Xt$ and whose restriction to $\Ut:=\Xt\stm\Dt$ is a variation of Hodge structure. Let $j_{\Ut}\,{:}\,\Ut\into\Xt$, $i_{\Dt}\,{:}\,\Dt\into\Xt$ be natural inclusions. In this case it may be natural to consider
\htt{2.3.1}{}
$$\aligned\M_{\Dt}:={}&\bl(\Gr^W_{w+1}\!H^1i_{\Dt}^!\M_{\Xt}\br)(1)\\={}&\bl(\Gr^W_{w+1}\!\bl((j_{\Ut})_*j_{\Ut}^*\M_{\Xt}\big/\M_{\Xt}\br)\br)(1).\endaligned
\leqno(2.3.1)$$
This is endowed with
$$\gamma\in{\rm Ext}^1\bl(\M_{\Dt}(-1),\Mt_{\Xt}\br),$$
using the short exact sequence
\htt{2.3.2}{}
$$0\to\Mt_{\Xt}\to(j_{\Ut})_*j_{\Ut}^*\M_{\Xt}\to H^1i_{\Dt}^!\M_{\Xt}\to 0,
\leqno(2.3.2)$$
where the direct image by $i_{\Dt}$ is omitted to simplify the notation. Its dual $\rho$ can be defined by using the self-duality isomorphisms of $\M_{\Xt}$, $\M_{\Dt}$ coming from their polarizations, where one can take the multiplicities of irreducible components of $\Dt$ into account. It does not seem, however, easy to prove that their compositions coincide with the action of the cohomology class of the divisor $\Dt$ (with multiplicities of irreducible components).
\sk
The examples below show that the coincidence of the action of the cohomology class of the divisor on $\ft_*\M_{\Xt}$ with the composition $\gamma\ssc\rho$ does {\it not\1} necessarily imply the coincidence of the action on $\ft_*\M_{\Dt}$ with $\rho\ssc\gamma$. (It is then quite nontrivial whether these two conditions are compatible with each other.)
\par\htt{E2.3a}\msn
{\bf Example~2.3a.} Let $X$ be a smooth projective variety of dimension $n\gess 3$ with $D\sst X$ a {\it smooth\1} divisor. Let $\M_X$ be a pure Hodge module of weight $n{+}1$ with strict support $X$ and whose restriction to $U\,{:=}\,X\stm D$ is a variation of Hodge structure of level~1, weight~1, and rank~2. (For instance, $X\eq Y{\times}C'$ with $\dim Y\eq n{-}1$, $\dim C'\eq 1$, and $\M_X$ is the (shifted) pullback of a pure Hodge module $\M_{C'}$ with strict support $C'$ which is associated with an elliptic surface over $C'$.)
\sk
Assume the local monodromy at $q\ins D$ is unipotent and non-semisimple. Let $\pi\,{:}\,\Xt\tos X$ be the blow-up at $q$ with $E\,(\cong\PP^{n-1})$ the exceptional divisor. Let $\Dt$ be the total transform of $D$, and $Z\,(\cong\PP^{n-2})$ be the intersection of $E$ with the proper transform of $D$. Set $\Ut\,{:=}\,\Xt\stm\Dt$ with $i_E\,{:}\,E\into\Xt$, $j_{\Ut}\,{:}\,\Ut\into\Xt$ natural inclusions. Let $\M_{\Xt}$ be the pure Hodge module of weight $n{+}1$ with strict support $\Xt$ and whose restriction to $\Ut\eq U$ is the given variation of Hodge structure. We have the following.
\par\htt{P2.3a}\msn
{\bf Proposition~2.3a.} {\it There is an isomorphism
\htt{2.3.3}{}
$$\M_{\Dt}(-1):=\Gr^W_{n+2}\bl((j_{\Ut})_*j_{\Ut}^*\M_{\Xt}/\M_{\Xt}\br)\cong\M_Z(-2),
\leqno(2.3.3)$$
where $\M_Z$ is a pure Hodge module of weight $n{-}2$ with strict support $Z$ and whose underlying $\R$-complex $K_Z$ is isomorphic to $\R_Z[n{-}2]$ {\rm(}since $Z$ is simply connected$\1)$.}
\msn
{\it Proof.} This is shown by using an approximative description of mixed Hodge modules of normal crossing type as in \cite[(3.18.5)]{mhm} (see also \cite{rh} for the underlying $\D$-modules). Indeed, cutting $\Xt$ by a (local) transversal slice $X'\sst\Xt$ to $Z$, we may replace $\Dt$ by a divisor $D'$ with normal crossings on a surface $X'$, and $\M_{\Xt}$ by the restriction $\M_{X'}$ of $\M_{\Xt}$ to $X'$ (up to a shift). The latter can be described approximatively by using the functors
$$\psi_x\psi_y,\q\psi_x\varphi_y,\q\varphi_x\psi_y,\q\varphi_x\varphi_y,$$
where $x,y$ are local coordinates of $X'$ with $D'\eq\{xy\eq0\}$ locally. Note that the stalk at the origin (that is, $H^{\ssb}i_0^*$ with $i_0\,{:}\,\{0\}\into X'$ a natural inclusion) of a mixed Hodge module of normal crossing type can be described by using the associated single complex of the double complex
\htt{2.3.4}{}
$$\begin{array}{cccccccc}\psi_x\psi_y&\to&\psi_x\varphi_y,\\ \downarrow&&\downarrow\\ \varphi_x\psi_y,&\to&\varphi_x\varphi_y,\end{array}
\leqno(2.3.4)$$
where the vertical and horizontal morphisms are respectively given by ${\rm can}_x$ and ${\rm can}_y$ (which correspond up to a sign to the actions of $\dd_x$ and $\dd_y$ on the underlying $\D$-module), see for instance \cite[Cor.\,2.24]{mhp}.
\sk
Applying the above four functors to $\M_{X'}$, we get respectively the $\R$-Hodge structures
\htt{2.3.5}{}
$$\bl[\R,\R(\mis 1)\br],\q\R(\mis 1),\q\R(\mis 1),\q0,
\leqno(2.3.5)$$
where $\bl[\R,\R(\mis 1)\br]$ denotes an extension of $\R(\mis 1)$ by $\R$, see also \cite{CKS2}, etc.
\sk
For the localization along $D'$ (that is, replacing $\M_{X'}$ with $(j_{X'\setminus D'})_*j_{X'\setminus D'}^*\M_{X'}$, where $j_{X'\setminus D'}\,{:}\,X'\stm D'\into X'$ is the inclusion), we obtain
\htt{2.3.6}{}
$$\bl[\R,\R(\mis 1)\br],\,\,\bl[\R(\mis 1),\R(\mis 2)\br],\,\,\bl[\R(\mis 1),\R(\mis 2)\br],\,\,\bl[\R(\mis 2),\R(\mis 3)\br],
\leqno(2.3.6)$$
since ${\rm var}_x$, ${\rm var}_y$ are bijective. (Note that ${\rm var}_x$, ${\rm var}_y$ correspond to the actions of $x,y$ on the underlying $\D$-module up to a sign.) Their quotient $\Cc$ corresponds to
\htt{2.3.7}{}
$$0,\q\R(\mis 2),\q\R(\mis 2),\q\bl[\R(\mis 2),\R(\mis 3)\br].
\leqno(2.3.7)$$
Here only $\R(\mis 2)$ in $\varphi_x\varphi_y\Cc$ corresponds to a subquotient of $\C$ with weight 4. (Note that the weight is shifted by 1 under the functors $\psi_x,\psi_y$, but it is unchanged by $\varphi_x,\varphi_y$.) So the isomorphism (\hl{2.3.3}{2.3.3}) follows, since the weight is shifted by $n{-}2\,({=}\,{\rm codim}_{\Xt}X')$ when we restrict to $X'\sst\Xt$. This finishes the proof of Proposition~\hl{P2.3a}{2.3a}.
\sk
The action of the cohomology class of $\Dt$ (with general multiplicities) on the cohomology of $\M_Z$ does not seem to vanish. (Indeed, the restriction of $\OO_{\Xt}(E)$ to $E\cong\PP^{n-1}$ is isomorphic to $\OO_{\PP^{n-1}}(-1)$.) We have, however, the following.
\par\htt{P2.3b}\msn
{\bf Proposition~2.3b.} {\it The direct image of the morphism
$$\gamma:\M_Z(-2)\to\M_{\Xt}[1]$$
by $\pi\,{:}\,\Xt\tos X$ vanishes.}
\msn
{\it Proof.} Taking the direct image of $\gamma$ by $\pi\,{:}\,\Xt\tos X$, we get the morphism
\htt{2.3.8}{}
$$\pi_*\gamma:\mopl_{i=0}^{n-2}\,\R_{\{q\}}(-i{-}2)[n{-}2{-}2i]\to\pi_*\M_{\Xt}[1]=\M_X[1],
\leqno(2.3.8)$$
since $K_Z\,{\cong}\,\R_Z[n{-}2]$. Here $\R_{\{q\}}$ denotes the Hodge module of weight 0 supported at $q$ and corresponding to the trivial Hodge structure $\R$, see Lemma~\hl{L2.3}{2.3} below for the last isomorphism of (\hl{2.3.8}{2.3.8}).
\sk
The morphism (\hl{2.3.8}{2.3.8}) must factor through $i_q^!\M_X$ with $i_q\,{:}\,\{q\}\into X$ the inclusion (using the adjunction relation), but
\htt{2.3.9}{}
$$H^ji_q^!\M_X=0\q(j\ne n),
\leqno(2.3.9)$$
since $H^ji_q^*\M_X=0\,\,\,(j\,{\ne}\,{-}n)$ (using the self-duality of $\M_X$). This implies the vanishing of (\hl{2.3.8}{2.3.8}), and we get Proposition~\hl{P2.3b}{2.3b}.
\sk
It remains to show the following.
\par\htt{L2.3}\msn
{\bf Lemma~2.3.} {\it The last isomorphism of {\rm(\hl{2.3.8}{2.3.8})} holds.}
\msn
{\it Proof.} We have the Leray spectral sequence
\htt{2.3.10}{}
$$E_2^{p,q}=H^p(E,\Hc^qK_{\Xt}|_E)\Longrightarrow H^{p+q}(E,K_{\Xt}|_E),
\leqno(2.3.10)$$
with $K_{\Xt}$ the underlying $\R$-complex of $\M_{\Xt}$. For $x\in E$, there are isomorphisms
\htt{2.3.11}{}
$$\Hc^q(K_{\Xt})_x\cong\begin{cases}\R&\h{if $\,q\eq{-}n$},\\ \R(-1)&\h{if $\,q\eq{-}n{+}1,\,x\ins Z$},\\0&\h{otherwise},\end{cases}
\leqno(2.3.11)$$
using (\hl{2.3.4}{2.3.4}--\hl{2.3.5}{5}). Since $E,Z$ are simply connected, these imply that
\htt{2.3.12}{}
$$E_2^{p,q}\cong\begin{cases}\R\bl(-\tfrac{p}{2}\br)&\h{if $\,q={-}n,\,\,p\in[0,2n{-}2]\cap 2\1\Z$},\\ \R\bl(-\tfrac{p}{2}{-}1\br)&\h{if $\,q={-}n{+}1,\,\,p\in[0,2n{-}4]\cap 2\1\Z$},\\0&\h{otherwise},\end{cases}
\leqno(2.3.12)$$
in particular, the spectral sequence degenerates at $E_3$. The picture of the spectral sequence is as below.
$$\setlength{\unitlength}{1cm}
\begin{picture}(11.7,1.5)
\put(.6,0){\h{$\q\R$}}
\put(.6,1){\h{$\R(-1)$}}
\put(2.1,.75){\vector(2,-1){.5}}
\put(2.1,.42){{$\scriptstyle{\rm d}_2$}}
\put(3,0){\h{$\R(-1)$}}
\put(3,1){\h{$\R(-2)$}}
\put(4.5,.75){\vector(2,-1){.5}}
\put(4.5,.42){{$\scriptstyle{\rm d}_2$}}
\put(5.4,0){\h{$\R(-2)$}}
\put(5.5,1.02){$\cdots$}
\put(7.5,.02){$\cdots$}
\put(6.9,1){\h{$\R(1{-}n)$}}
\put(8.7,.75){\vector(2,-1){.5}}
\put(8.7,.42){{$\scriptstyle{\rm d}_2$}}
\put(9.5,0){\h{$\R(1{-}n)$}}
\end{picture}$$
\sk
By the decomposition theorem, the cohomology groups $H^j(E,K_{\Xt}|_E)$ must have {\it pure weight\1} $j{+}n{+}1$ for $j\ne -n$. Indeed, $\M_{\Xt}$ has pure weight $n{+}1$, $\pi$ induces an isomorphism over $X\stm\{q\}$, and $\Hc^jK_X\eq0$ for $j\,{\ne}\,{-}n$ (because $D$ is smooth), where the ``classical $t$-structure" is used, see \cite[4.6]{mhm}. Then the differential $\ddd_2\,{:}\,E_2^{p,q}\tos E_2^{p+2,q-1}$ must be surjective (hence bijective) for $(p,q)\eq(2i,-n{+}1)$ ($i\ins[0,n{-}2]$ using the $E_3$-degeneration. We then get the vanishing of $E_3^{p,q}$ except for $(p,q)\eq(0,-n)$, and the last isomorphism of (\hl{2.3.8}{2.3.8}) follows. This completes the proofs of Lemma~\hl{L2.3}{2.3} and Proposition~\hl{P2.3b}{2.3b}.
\par\htt{E2.3b}\msn
{\bf Example~2.3b.} Let $X$ be a smooth projective variety of dimension $n\gess 3$ with $D\sst X$ a divisor with simple normal crossings such that the singular locus $\Si$ of $D$ is smooth. Let $\M_X$ be a pure Hodge module of weight $n{+}1$ with strict support $X$ and whose restriction to $U\,{:=}\,X\stm D$ is a variation of Hodge structure of level~1, weight~1, and rank~2. (For instance, $X$ is a blow-up of $Y{\times}C'$ in Example~\hl{E2.3a}{2.3a}.)
\sk
Assume the local monodromies around any smooth points of $D$ are unipotent and non-semisimple. Let $\pi\,{:}\,\Xt\tos X$ be the blow-up at a point $q\ins\Si$ with $E\,(\cong\PP^{n-1})$ the exceptional divisor. Let $\Dt$ be the total transform of $D$ with $Z_0$ the strict transform of $\Si$. Let $Z_1\cup Z_2$ (with $Z_k\cong\PP^{n-2})$ be the intersection of $E$ with the proper transform of $D$ (which is locally the union of $D_1$ and $D_2$). Set $\Ut\,{:=}\,\Xt\stm\Dt$ with $i_E\,{:}\,E\into\Xt$, $j_{\Ut}\,{:}\,\Ut\into\Xt$ natural inclusions. Let $\M_{\Xt}$ be the pure Hodge module of weight $n{+}1$ with strict support $\Xt$ and whose restriction to $\Ut\eq U$ is the given variation of Hodge structure. In this case we have
\htt{2.3.13}{}
$$\M_{\Dt}\cong\mopl_{k=0}^2\,\M_{Z_k}(-1),
\leqno(2.3.13)$$
where the $\M_{Z_k}$ are pure Hodge modules of weight $n{-}2$ with strict support $Z_k$ and their underlying $\R$-complexes $K_{Z_k}$ are isomorphic to $\R_{Z_k}[n{-}2]$ (on a neighborhood of $E$).
\sk
Set $Z^o_k:=Z_k\stm\Si$ ($k\eq1,2$). For $x\ins E$, we have
\htt{2.3.14}{}
$$\Hc^j(K_{\Xt})_x\cong\begin{cases}\R&\h{if $\,j\eq{-}n$},\\ \R(-1)&\h{if $\,j\eq{-}n{+}1,\,x\ins Z^o_1\cup Z^o_2$},\\\R^2(-1)&\h{if $\,j\eq{-}n{+}1,\,x\ins Z_0\cap E$},\\0&\h{otherwise},\end{cases}
\leqno(2.3.14)$$
hence
\htt{2.3.15}{}
$$\Hc^j(K_{\Xt})|_E\cong\begin{cases}\R_E&\h{if $\,j\eq{-}n$},\\ \R_{Z_1}(\mis 1)\oplus\R_{Z_2}(\mis 1)&\h{if $\,j\eq{-}n{+}1$},\\0&\h{otherwise}.\end{cases}
\leqno(2.3.15)$$
Indeed, applying the four functors in (\hl{2.3.4}{2.3.4}) to $K_{\Xt}|_E$ restricted to a transversal slice, we get
$$\R,\q\R(\mis 1),\q\R(\mis 1),\q0.$$
Note that the restriction to $E$ can be expressed by the mapping cone of ${\rm can}_z\,{:}\,\psi_{z,1}\tos\varphi_{z,1}$ if $E\eq\{z\eq0\}$ locally. Here ${\rm can}_z$ corresponds up to a sign to the action of $\dd_z$ on the underlying $\D$-module.
\sk
By an argument similar to Example~\hl{E2.3a}{2.3a}, we then conclude that
\htt{2.3.16}{}
$$\aligned\pi_*\M_{\Xt}&\cong\M_X\oplus\mopl_{i=0}^{n-3}\,\R_{\{q\}}(-i{-}2)[n{-}3{-}2i],\\
\pi_*\M_{\Dt}&\cong\M_{\Si}(-1)\oplus\mopl_{i=0}^{n-4}\,\R_{\{q\}}(-i{-}2)[n{-}4{-}2i]\\&\q\oplus\mopl_{i=0}^{n-2}\,\R^2_{\{q\}}(-i{-}1)[n{-}2{-}2i],\endaligned
\leqno(2.3.16)$$
where $\R^2_{\{q\}}\eq\R_{\{q\}}\oplus\R_{\{q\}}$. (Note that $Z_0$ is a one-point blow-up of $\Si$.) This implies that, in order to get the coincidence of the action of the cohomology class $\ell'$ of the divisor with the compositions of $\rho,\gamma$, one may have to {\it delete\1} some direct factor (compatible with the self-duality and also with the action of $\ell'$) from the direct image of $\M_{\Dt}$, considering the difference in the {\it rank\1} of the action of $\ell'$ on $\pi_*\M_{\Xt}$ and $\pi_*\M_{\Dt}$ when $n$ is not very small. This {\it cannot\1} be done on $\Xt$ before taking the direct image. Here one cannot delete simply the last term on the right-hand side of the second isomorphism of (\hl{2.3.16}{2.3.16}), since one might then get a $\M_{\Dt}$-decomposable relation, see Remark~\hl{R2.2}{2.2}. (This is closely related to the image of the morphism from $H^{\ssb}\pi_*\Gr^W_{w+2}H^1i_{\Dt}^!\M_{\Xt}$. However, this image is not a direct factor compatible with the self-duality and the action of $\ell'$, and it becomes too big if we consider the smallest direct factor containing it and compatible with the self-duality and the action of $\ell'$.) One may have to take also $\Gr^L_{\ssb}$ in the case $\Xt$ is replaced by its product with a smooth curve (that is, if the center of blow-up becomes a curve).
\par\htt{R2.3a}\msn
{\bf Remark~2.3a.} One may have to use a point-center blow-up when one considers for instance the subspace $X{\times}\{q\}\cup\{q\}{\times}X\subset X{\times}X$ if desingularization is obtained by repeating blow-ups with smooth centers as in \cite{Wl} without using normalization. (Note that normalization is not necessarily compatible with closed embeddings into complex manifolds in general.) There is an example having a self-intersection if we consider for instance the quotient of
$$X{\times}\{q\}\cup\{q\}{\times}X\cup X{\times}\{q'\}\cup\{q'\}{\times}X\subset X{\times}X$$
by an automorphism subgroup generated by $(x,y)\mapsto(y,\iota(x))$ (deleting the singular points), where $\iota$ is an involution of $X$ such that the fixed point set of $\iota$ is discrete and $\iota(D)\eq D$, $\iota(q)\eq q'\,{\ne}\,q$. (One can replace $C'$ with a double cover, and assume $Y$ is a complex torus, for instance.)
\par\htt{R2.3b}\msn
{\bf Remark~2.3b.} There seems to be a certain problem even in the constant Hodge module case if $\Dt$ is reducible. Consider for instance the case $X\eq\PP^3$, $D\eq\PP^2\sst X$, $\Xt$ is the blow-up of $X$ along $C\,{:=}\,\PP^1\sst D$ with exceptional divisor $E\eq\PP^1{\times}\PP^1$. Let $\Dt$ be the total transform of $D$, which is the union of $E$ and the proper transform $D'\,({\cong}\,D)$ intersecting transversally along $C\sst D'$, since $C$ has codimension~1 in $D$. The cohomology groups of degree $k$ of $\Xt,E,D',C$ are direct sums of $\R(-k/2)$ if $k$ is even, and vanish otherwise. Their dimensions for even degrees are respectively as follows:
$$(1,2,2,1),\q(1,2,1),\q(1,1,1),\q(1,1).$$
The smallest graded vector subspace of $H^{\ssb}(E)\,{\oplus}\,H^{\ssb}(D')$ compatible with the self-duality and the action of $\ell'$ (which vanishes on $H^{\ssb}(D')$) and containing the image of $H^{\ssb-2}(C)(-1)$ under the Gysin morphism seems to be the whole space (since the intersection form of the middle cohomology of $E$ is hyperbolic).

\end{document}